\documentclass[a4paper,12pt]{article}

\usepackage{amsmath}
\newcommand{\modelo}[1]{#1}
\newcommand{\variedad}[1]{\mathcal{#1}}

\newcommand{\A}{\modelo{A}}

\newcommand{\V}{\variedad{V}}

\usepackage{amsthm}
\usepackage{verbatim}
\newlength{\ancho}
\textwidth=\ancho
\newcommand{\vu}{\vec u}
\newcommand{\vv}{\vec v}
\newcommand{\vw}{\vec w}

\newcommand{\<}{\langle}
\renewcommand{\>}{\rangle}
\newcommand{\ent}{\Rightarrow}
\newcommand{\tne}{\Leftarrow}

\renewcommand{\phi}{\varphi}

\newtheorem{theorem}{Theorem}

\newtheorem{corollary}[theorem]{Corollary}
\newtheorem{definition}[theorem]{Definition}

\begin{document}
\title{Existentially Definable Factor Congruences}
\author{Pedro S\'{a}nchez Terraf\thanks{Supported by CONICET\newline \emph{2000 Mathematics Subject Classification:} Primary
  08B05, Secondary 03C40.
}}
\date{}
\maketitle
\begin{abstract}
A variety $\V$  has \emph{definable factor congruences} if and only if
factor congruences can be defined by a first-order formula $\Phi$ having
\emph{central elements} as parameters. We prove that if $\Phi$ 
can be chosen to be existential, factor congruences in every
algebra of $\V$ are compact.
\end{abstract}

We study factor congruences in order to understand direct product
representations in varieties. It is known that in rings with identity
and bounded lattices, factor congruences are characterized,
respectively, by central idempotent elements and neutral complemented
elements. D. Vaggione~\cite{va0} generalized these concepts to a broader
context. A \emph{variety with $\vec{0}$ \& $\vec{1}$} is a variety $\V$ in which
there exist unary terms $0_{1}(w),\dots,0_{l}(w),$ $1_{1}(w),\dots,1_{l}(w)$ such that
\begin{equation*}
\mathcal{V}\models \vec 0(w)=\vec 1(w)\rightarrow x=y,
\end{equation*}
where $w$, $x$ and $y$ are distinct variables,
$\vec{0}=(0_{1},\dots,0_{l})$ and $\vec{1}=(1_{1},\dots,1_{l}).$ If
$\lambda \in A\in \mathcal{V}$, we say that $\vec{e}\in A^{l}$ is a \emph{$\lambda$-central element} of $A$
if there exists an isomorphism $A\rightarrow A_{1}\times A_{2}$ such
that 
\[\lambda \mapsto \<\lambda_1,\lambda_2\>,\]
\[\vec{e}\mapsto [ \vec{0}(\lambda_1),\vec{1}(\lambda_2)].\]
where we write $[\vec a, \vec b]$ in place
of $ (\<a_{1},b_{1}\>,\dots,\<a_{l},b_{l}\>) \in (A \times B)^l$ for $\vec a \in A^l$ and $\vec b \in B^l$.
It is clear from the above definitions that if the language of $\V$
has a constant symbol $c$, the terms $\vec{0}$ and  $\vec{1}$ can be chosen
closed, and we can define a \emph{central element} of $A$ to be just a
$c^A$-central element. We will work heretofore under this assumption.

In~\cite{DFC}, Vaggione and the author
introduced the following concept:
\begin{definition}
 $\V$ has \emph{Definable Factor Congruences (DFC)} iff there exists a first order formula $\Phi(x,y,\vec z)$ in the
  language of $\V$ such that for all $A, B\in \V$, and $a,c\in A$,
  $b,d\in B$,
\begin{equation}\label{eq:DFC1}
A\times B \models \Phi\bigl(\<a,b\>, \<c,d\>, [\vec 0, \vec 1]\bigr)
  \quad \text{ if and only if } \quad a=c.
\end{equation}
\end{definition}
Varieties with $\vec{0}$ \& $\vec{1}$ are
obviously \emph{semidegenerate} (no non-trivial algebra in the variety has a
trivial subalgebra). In~\cite{DFC} it is  proved that DFC is equivalent to Boolean factor
congruences \cite{7,4,ChaJonTar} in semidegenerate varieties. 

It is then natural to ask if, for a variety $\V$ with DFC, the
quantifier complexity of $\Phi$ is reflected in the structure of
$\V$. One early work \cite{9} showed that if factor congruences in
$\V$ are compact, then $\V$ has DFC and $\Phi$ can be chosen to be a
positive existential formula. A partial converse was proved in
\cite{DFC}:  if a positive formula $\Phi$ witnesses DFC  for $\V$,
$\V$ has compact factor congruences. Also, counterexamples were
constructed showing that for universal $\Phi$ this may not be the case.
  
In this note we prove:
\begin{theorem}\label{th:main}
Let $\V$ be a variety with   $\vec{0}$ \& $\vec{1}$. Suppose there exists an existential formula $\Phi$ 
that satisfies~(\ref{eq:DFC1}).  Then we may replace $\Phi$ by a
positive formula.
\end{theorem}
By using the results cited from \cite{DFC} and \cite{9} we obtain
\begin{corollary}
Let $\V$ be any variety. The following are equivalent:
\begin{enumerate}
\item There exist  unary terms $\vec 0(w) = 0_{1}(w),\dots,0_{l}(w),$
  $\vec 1(w) = 1_{1}(w),\dots,1_{l}(w)$ and 
 an existential first-order formula $\Phi(x,y,\vec z)$ in the
  language of $\V$ such that for all $A, B\in \V$, and $a,c,e\in A$,
  $b,d,f\in B$,
  \[A\times B \models \Phi\bigl(\<a,b\>, \<c,d\>, [\vec 0(e), \vec 1(f)]\bigr) 
  \quad \text{ if and only if } \quad a=c.\] 
\item $\V$ has compact factor congruences.
\end{enumerate}
\end{corollary}
\begin{proof}[Proof of Theorem~\ref{th:main}]
We will only consider the case $l=1$, so we have two closed terms 0 and 1
that satisfy
\begin{equation}\label{eq:DFC}
A\times B \models \Phi\bigl(\<a,b\>, \<c,d\>, \<0, 1\>\bigr)
  \quad \text{ if and only if } \quad a=c.
\end{equation}
The general case is straightforward.

We will write $F(x_1,\dots,x_n)$ for the algebra  freely generated by
$\{x_1,\dots,x_n\}$ in $\V$. Assume
\begin{equation*}\label{eq:1}
\Phi (x,y,z):= \exists \vw \bigvee_{i} \bigwedge_{j}
\phi_{ij}(x,y,z,\vw)
\end{equation*}
where $\phi_{ij}(x,y,z,\vu)$ is atomic or negated atomic and let  $\Lambda_i = \{ j : \phi_{ij} \text{ is atomic}\}$.
Taking $A:=F(x)$ and $B:= F(x,y)$ in~(\ref{eq:DFC}) we obtain:
\[F(x) \times F(x,y)   \models  \exists \vw \Phi(\<x,x\>, \<x,y\>, \<0,1\>
,\vw)\]
Hence there exists $ [\vu(x),\vv(x,y)]$ in $F(x) \times F(x,y)$
and $k$ such that
\begin{equation}\label{eq:4}
F(x) \times F(x,y)   \models  \bigwedge_{j} \phi_{kj}(\<x,x\>, \<x,y\>,
\<0,1\>,  [\vu(x),\vv(x,y)])
\end{equation}
Using  preservation by homomorphic images,  we obtain
\begin{equation}\label{eq:2}
\begin{split}
\V &\models  \bigwedge_{j\in \Lambda_k} \phi_{kj}(x,x,0,  \vu(x))\\
\V &\models  \bigwedge_{j\in \Lambda_k} \phi_{kj}(x,y,1,  \vv(x,y))
\end{split}
\end{equation}
Now we will prove that for this $k$, the positive formula
\begin{equation*}\label{eq:3}
\Phi'(x,y,z):= \exists \vw   \bigwedge_{j\in \Lambda_k} \phi_{kj}(x,y,z,\vw)
\end{equation*}
satisfies~(\ref{eq:DFC}).

\noindent ($\tne$) Take  $a\in\A$, $b,c\in B$, $A,B\in\V$. Using~(\ref{eq:2})
and preservation by direct products, 
\[A\times B \models  \bigwedge_{j\in \Lambda_k} \phi_{kj}(\<a,b\>, \<a,c\>,
\<0,1\>,  [\vu(a),\vv(b,c)] ).\]
Hence
\[A\times B \models \exists \vw \bigwedge_{j\in \Lambda_k} \phi_{kj}(\<a,b\>, \<a,c\>,
\<0,1\>,  \vw ),\]
and by definition,
\[A\times B \models  \Phi'(\<a,b\>, \<a,c\>,\<0,1\>) \]
\medskip

\noindent ($\ent$) Now suppose $A\times B \models  \Phi'(\<a,b\>,
\<c,d\>,\<0,1\>)$. By 
preservation by homomorphic images, we have $A \models \Phi'(a,c,0)$; take
$\vw$ such that
\begin{equation}\label{eq:5}
A\models  \bigwedge_{j\in \Lambda_k} \phi_{kj}(a,c,0,\vw).
\end{equation}

Considering (recall~(\ref{eq:4}))
\[ F(x)\times F(x,y) \models \bigwedge_{j\not\in \Lambda_k} \phi_{kj}(\<x,x\>, \<x,y\>,
\<0,1\>,  [\vu(x),\vv(x,y)] ) \]
we obtain
\[A\times \bigl(F(x)\times  F(x,y)\bigr) \models \bigwedge_{j\not\in \Lambda_k} \phi_{kj}\bigl(\bigl\<a,\<x,x\>\bigr\>, \bigl\<c,\<x,y\>\bigr\>,
\bigl\<0,\<0,1\>\bigr\>,  \bigl[\vw, [\vu(x),\vv(x,y)] \bigr] \bigr), \]
and hence, using~(\ref{eq:4}) and~(\ref{eq:5}),
\[A\times \bigl(F(x)\times  F(x,y)\bigr) \models \bigwedge_{j} \phi_{kj}\bigl(\bigl\<a,\<x,x\>\bigr\>, \bigl\<c,\<x,y\>\bigr\>,
\bigl\<0,\<0,1\>\bigr\>,  \bigl[\vw, [\vu(x),\vv(x,y)] \bigr] \bigr). \]
Equivalently, by the obvious isomorphism
\[\bigl(A\times F(x)\bigr)\times  F(x,y) \models \bigwedge_{j}
\phi_{kj}\bigl(\bigl\<\<a,x\>,x\bigr\>,
\bigl\<\<c,x\>,y\bigr\>,\bigl\<\<0,0\>,1\bigr\>,  \bigl[[\vw,
  \vu(x)],\vv(x,y) \bigr] \bigr). \]
This yields, taking, $\mathbf{e}:= \<a,x\>$, $\mathbf{f} := \<c,x\>$,
$\mathbf{0} := \<0,0\> = 0^{A\times F(x)}$, 
\[\bigl(A\times F(x)\bigr)\times  F(x,y) \models \Phi( \<\mathbf{e},x\>,
\<\mathbf{f},x\> , \<\mathbf{0} ,1\>)\]
and by~(\ref{eq:DFC}),
\[A\times F(x) \models \mathbf{e} =\mathbf{f}.\]
I.e., $\<a,x\> = \<c,x\>$. We may conclude $a=c$, as desired.
\end{proof}

\begin{quote}
CIEM --- Facultad de Matem\'atica, Astronom\'{\i}a y F\'{\i}sica 
(Fa.M.A.F.) 

Universidad Nacional de C\'ordoba --- Ciudad Universitaria

C\'ordoba 5000. Argentina.

\texttt{sterraf@famaf.unc.edu.ar}
\end{quote}
\end{document}